\documentclass[a4paper,12pt]{amsart}
\usepackage[pdftex]{color,graphicx}
\usepackage{amssymb}
\usepackage{amsfonts}
\usepackage{esint}
\usepackage{relsize}
\usepackage{amsmath}
\usepackage{amsrefs}
\usepackage{epsfig}
\usepackage{amsthm}
\usepackage{color}
\usepackage[]{epsfig}
\usepackage[]{pstricks}
\usepackage[utf8]{inputenc}
\usepackage{multicol}
\usepackage{accents}
\newtheorem{theorem}{Theorem}[section]

\newtheorem{corollary}{Corollary}[section]

\theoremstyle{definition}
\newtheorem{definition}{Definition}[section]

\theoremstyle{remark}
\newtheorem{remark}{Remark}[section]
\newcommand{\R}{\mathbb{R}}
\newcommand{\h}{\mathbb{H}}
\newcommand{\s}{\mathbb{S}}

\newcommand{\CC}{\mathbb{C}}

\newcommand{\ria}{\rightarrow}

\newcommand{\om}{\omega}
\newcommand{\n}{\nabla}

\newcommand{\ran}{\rangle}
\newcommand{\lan}{\langle}
\newcommand{\ve}{\varepsilon}
\newcommand{\vp}{\varphi}

\DeclareMathOperator{\rad}{rad}

\DeclareMathOperator{\tr}{trace}

\numberwithin{equation}{section}
\title{Hopf type theorems in Riemannian manifolds}

\author{Hil\'ario Alencar, Greg\'orio Silva Neto \& Detang Zhou}

\dedicatory{Dedicated to Professor Renato Tribuzy on the occasion of his 75th birthday}

\date{November 08, 2021}

\address{Instituto de Matemática, Universidade Federal de Alagoas, Macei\'o, AL, 57072-900, Brazil}

\email{hilario@mat.ufal.br}

\address{Instituto de Matemática, Universidade Federal de Alagoas, Macei\'o, AL, 57072-900, Brazil}

\email{gregorio@im.ufal.br}

\address{Instituto de Matemática e Estatística, Universidade Federal Fluminense, Niterói, RJ, 24020, Brazil}

\email{zhoud@id.uff.br}

\begin{document}
%\subjclass[2020]{Primary  53E10; Secondary 53C42; 53C21; 35C06; 35A15; 35A23; 35J15; 35J60}

 
\footnotetext{Hil\'ario Alencar, Greg\'orio Silva Neto and Detang Zhou were partially supported by the National Council for Scientific and Technological Development - CNPq of Brazil.}% Manuel Cruz was partially supported by Coordena\c c\~ao de Aperfei\c coamento de Pessoal de N\'ivel Superior - Brasil (CAPES) - Finance Code 001.}

\begin{abstract}
In 1951, H. Hopf proved that the only surfaces, homeomorphic to the sphere, with constant mean curvature in the Euclidean space are the round (geometrical) spheres. In this paper we survey some contributions of Renato Tribuzy to generalize Hopf's result as well as some recent results of the authors using these techniques for shrinking solitons of curvature flows and for surfaces in three-dimensional warped product manifolds, specially the de Sitter-Schwarzschild and Reissner-Nordstrom manifolds.
\end{abstract}
\maketitle

\section{Introduction and Tribuzy's contributions}

In 1951, H. Hopf, see \cite{hopf} and \cite{hopf-1000}, proved that the only surfaces with constant mean curvature in $\R^3$, homeomorphic to the sphere, are the round spheres. Hopf gave two different proofs of this result. In the first proof, one considers the second fundamental form $II$ in isothermal parameters and takes the $(2,0)$-component of $II$, i.e., $II^{(2,0)}=(1/2)P dz^2.$ It can be shown that the complex function $P$ vanishes precisely at the umbilical points of $\Sigma$ and it is holomorphic if and only if the mean curvature of $\Sigma$ is constant. It is also seen that the quadratic form $II^{(2,0)}$ does not depend on the parameter $z$; hence, it is globally defined on $\Sigma$. It is a known theorem on Riemann surfaces that if the genus $g$ of $\Sigma$ is zero, any holomorphic quadratic form vanishes identically. Then $P=0$, i.e., all points of $\Sigma$ are umbilic, and hence $\Sigma$ is a standard sphere. His second proof is based on the lines of curvature. The quadratic equation $\mbox{Im}(Pdz^2)=0$ determines two fields of directions (the principal directions), whose singularities are the zeroes of $P$. Since $P$ is holomorphic, if $z_0$ is a zero of $P$, either $P=0$ in a neighborhood $V$ of $z_0$ or 
\begin{equation}\label{eq-001}
P(z)=(z-z_0)^kh_k(z), \ z\in V, \ k\geq 1,
\end{equation}
where $h_k$ is a function of $z$ with  $h_k(z_0)\neq 0,$ see for example \cite{Rudin}, p. 208-209. This number $k$ is called the order of the zero. In particular, if $P$ is not identically zero in a neighborhood of $z_0,$ then $z_0$ is an isolated singularity of the field of directions with index $-k/2$. Since $\Sigma$ has genus zero, by the Poincar\'e index theorem, the sum of the indices of all singularities for any field of directions is two (hence positive). This lead us to a contradiction, and thus $P$ is identically zero.

In 1933, Carleman \cite{Carleman} was the first to observe that this property holds for non-analytic smooth functions which satisfies some first order partial differential equation. In fact, he proved that a solution $h:U\subset\mathbb{C}\rightarrow\mathbb{C}$ of 
\[
\frac{\partial h}{\partial\bar{z}} = ah+b\bar{h},
\]
does not admits a zero of infinite order except if $h=0.$ Here bars mean complex conjugate and $a,b$ are continuous complex functions. Notice that, if $a=b=0,$ then $h$ is holomorphic. Using these ideas, Hartman and Wintner, see \cite{H-W} and \cite{H-W-2}, and Chern, see \cite{chern}, proved their well known results on the classification of special Weingarten surfaces.

Following the ideas of Chern, in 1987 \cite{E-T-1} (see also \cite{E-T}), Eschenburg and Tribuzy proved the following result:

\begin{theorem}[Eschenburg-Tribuzy]\label{CR-ET}
Let $h:U\subset\mathbb{C}\to\mathbb{C}$ be a complex function defined in an open set $U$ of the complex plane. Assume that
\begin{equation}\label{CR}
\left|\frac{\partial h}{\partial \bar{z}}\right|\leq \vp(z)|h(z)|
\end{equation}
where $\vp$ is a $L^p,$ $p>2,$ non-negative real function. Assume further that $z=z_0\in U$ is a zero of $h.$ Then either $h\equiv0$ in a neighborhood $V\subset U$ of $z_0,$ or
\[
h(z)=(z-z_0)^kh_k(z),\ z\in V, \ k\geq 1,
\]
where $h_k(z)$ is a continuous function with $h_k(z_0)\neq 0.$
\end{theorem}

\begin{remark}
Eschenburg and Tribuzy called inequality \eqref{cauchy-E-T} the \emph{Cauchy-Riemann inequality}.
\end{remark}

By using this result, Eschenburg and Tribuzy extended the result of Hopf in the following way:

\begin{theorem}[Eschenburg-Tribuzy]\label{E-T-0}
Let $Q_c^3$ be a three-dimensional Riemannian manifold with constant sectional curvature $c\in\R.$ Let $X:\Sigma\ria Q_c^3$ be an immersed surface with mean curvature function $H.$ Assume that $\Sigma$ is homeomorphic to the sphere. If there exists a locally $L^p,$ $p>2,$ function $f:\Sigma\ria\R$ such that
\begin{equation}\label{cauchy-E-T}
|dH|\leq f\sqrt{H^2-K + c},
\end{equation}
where $K$ is the Gaussian curvature of $\Sigma,$ then $X(\Sigma)$ is totally umbilical.
\end{theorem}

The idea behind the proof is that the Hopf quadratic differential on an immersed sphere vanishes when the mean curvature is constant, and by definition of the Hopf differential its zeroes are umbilical points of the surface. Thus the immersed sphere is umbilical. When the mean curvature satisfies a Cauchy-Riemann inequality, then the zeros of the Hopf differential are all isolated and of negative index, or every point is a a zero of the Hopf form. This index is with respect to the two eigendirection lines at the points where the Hopf form is non zero. Since the sum of the indices of the zeros (when they are all isolated) is positive, the Hopf form must vanish on the surface. Thus each point is umbilical.

%\section{Generalization to other ambient spaces}

%In \cite{AR}, Abresch and Rosenberg found a quadratic form whose $(2,0)$-part is holomophic for constant mean curvature surfaces in $\mathbb{S}^2\times\R$ and $\mathbb{H}^2\times\R.$ 

When an immersed sphere in $M=\mathbb{S}^2\times\R$ or $M=\mathbb{H}^2\times\R$ has constant mean curvature, Abresch and Rosenberg, see \cite{AR}, found a holomorphic quadratic differential $Q$ on the surface generalizing the Hopf differential, that has the property that when it vanishes on the surface, then the sphere must be a rotationally invariant surface in $M$ (umbilical only when the surface is a slice, where the mean curvature is zero). This differential was found by calculating the second fundamental form of the rotational examples and making a simple modification of the Hopf form to find a quadratic form that vanishes on the rotational surfaces and that characterizes the rotational surfaces, i.e., when the form vanishes on an immersed surface, then the surface is indeed rotational.

To weaken the hypothesis of constant mean curvature, Alencar, do Carmo, and Tribuzy found in \cite{AdCT} a Cauchy-Riemann inequality involving the differential $Q$. This implies the zeros of $Q$ satisfy the same conditions as that of the Hopf form in the previous paragraph. Here the index of an isolated zero of $Q$ is with respect to the two eigendirections of $Q$ at a non zero point. Since the sum of the indices of the isolated zeros would be positive, $Q$ must vanish identically. Then one applies the theorem of Abresch-Rosenberg to conclude it is rotational. Namely, they obtained

%By using Theorem \ref{E-T-0}, the first author, do Carmo and Tribuzy \cite{AdCT} were able to prove the following

\begin{theorem}[Alencar-do Carmo-Tribuzy, 2007]\label{AdCT-2007}
Let $\Sigma$ be a compact immersed surface of genus zero in $\mathbb{H}^2\times\R$ or $\mathbb{S}^2\times\R.$ Assume that
\[
|dH|\leq \vp|Q^{(2,0)}|,
\]
where $|dH|$ is the norm of the differential $dH$ of the mean curvature $H$ of $\Sigma$, $Q^{(2,0)}$ is the $(2,0)$-part of the Abresch-Rosenberg quadratic differential, and $\vp$ is a continuous, non-negative real function. Then $Q^{(2,0)}$ is identically
zero, and $\Sigma$ is an embedded surface invariant by rotations in $\mathbb{H}^2\times\R$ or $\mathbb{S}^2\times\R.$
\end{theorem}

In \cite{AdCT2}, the first author, do Carmo and Tribuzy generalized the Abresch-Rosenberg quadratic differential to higher-codimensions and, using Theorem \ref{E-T-0}, the index argument explained earlier, and some arguments of reduction of codimension, as Theorem 4 of \cite{Yau}, they were able to prove

\begin{theorem}[Alencar-do Carmo-Tribuzy, 2010]
Let $\Sigma$ be a compact surface of genus zero and let $x:\Sigma\to E_c^n\times\R$, $n\geq 2,$ be an immersion of $\Sigma$ with parallel mean curvature, where $E_c^n$ is a space form of constant sectional curvature $c\in\R$. Then, one of the following assertions holds:
\begin{itemize}
\item[1)] $x(\Sigma)$ is a minimal surface of a totally umbilical hypersurface of $E_c^n$;
\item[2)] $x(\Sigma)$ is a standard sphere of a totally umbilical 3-dimensional submanifold of $E_c^n$;
\item[3)] $x(\Sigma)$ is a standard sphere of $E_c^3$;
\item[4)] $x(\Sigma)$ lies in $E_c^4\times\R\subset\R^6$ (possibly with the Lorentz metric), and there exists a plane $P$ such that $x(\Sigma)$ is invariant for rotations which fix its orthogonal complement. Furthermore, the level curves of the height
function $p\mapsto\lan x(p),\xi\ran$ are circles lying in planes parallel to $P$. Here, $\xi$ is the direction vector of the $\R$ component of $E_c^n\times\R.$
\end{itemize}
\end{theorem}

To conclude the list results with the contribution of Renato Tribuzy in this subject, we can also cite the result of the first author, do Carmo, Fernandez and Tribuzy \cite{AdCFT}. They generalized the Abresch-Rosenberg quadratic differential for three-dimensional simply connected homogeneous spaces with four dimensional isometry group $\mathbb{E}^3(\kappa,\tau).$ By using Theorem \ref{E-T-0} and an argument analogous to the proof of Theorem \ref{AdCT-2007}, they proved

\begin{theorem}[Alencar-do Carmo-Fernandez-Tribuzy, 2007]
Let $\Sigma$ be a compact surface of genus zero immersed in $\mathbb{E}^3(\kappa,\tau).$ Assume that
\[
|dH|\leq\vp|Q^{(2,0)}|,
\]
where $\vp$ is a non-negative real continuous function and $Q^{(2,0)}$ is the $(2,0)$ part of the generalized Abresch-Rosenberg quadratic differential in $\mathbb{E}^3(\kappa,\tau).$ Then $Q^{(2,0)}$ is identically zero and, by \cite{AR-2}, $\Sigma$ is a constant mean curvature surface invariant by rotations in $\mathbb{E}^3(\kappa,\tau).$
\end{theorem}

%\begin{remark}
%In this subject we can also cite Tribuzy's thesis where, among other results, he proved the following equivalence for quadratic forms
%\begin{lemma}
%Let $T$ be a symmetric bilinear form on a connected surface
%\end{lemma}
%\end{remark}

%Citar a tese do Renato Tribuzy em algum lugar, num remark.

\section{Weighted Riemannian manifolds}

In this section, we present a generalization of Theorem \ref{E-T-0} and its application to surfaces in $\mathbb{R}^n$ with weighted measure as we define later. These results are the main results the the paper \cite{ANZ-1} of the authors. 

An immersion $X:\Sigma\to\R^n$ of a two-dimensional surface $\Sigma$ is called a self-shrinker for the mean curvature flow if its mean curvature vector ${\bf H}$  satisfies the equation
\[
{\bf H}=-\frac{1}{2}X^\perp,
\]
where $X^\perp$ is the normal part of the position vector. Self-shrinkers are the self-similar solutions of the mean curvature flow and many efforts were made in the last decades in order to obtain examples of such surfaces and classify these surfaces under certain geometrical restrictions. In particular, there is a problem to classify the sphere as the only compact self-shrinker under some geometrical assumptions, as we can see, for example, in \cite{huisken}, \cite{C-M}, \cite{Cao-Li}, \cite{K-M}, \cite{brendle}, among others. In this spirit, the authors proved in \cite{ANZ-1} the following result:

\begin{theorem}[Theorem 1.2 of \cite{ANZ-1}]\label{Gauss-space-2}
Let $X:\Sigma\to\R^3$ be an immersed self-shrinker homeomorphic to the sphere. If there exists a non-negative locally $L^p$ function $\vp:\Sigma\to\R,\ p>2,$ and a locally integrable function $G:[0,\infty)\rightarrow[0,\infty)$ satisfying $\limsup_{t\rightarrow 0}G(t)/t<\infty,$ such that  
\begin{equation}\label{HG-2}
(\|X\|^2-4H^2)H^2\leq \vp^2 G(\|\Phi\|)^2,
\end{equation}
then $X(\Sigma)$ is a round sphere of radius $2$ and centered at the origin.

Here  $\|\Phi\|$ denotes the matrix norm of $\Phi=A-(H/2)I,$ where $A$ is the shape operator of the second fundamental form of $X,$ $H$ is its non-normalized mean curvature, and $I$ is the identity operator of $T\Sigma.$
\end{theorem}

Since the Hopf's differential is not holomorphic for self-shrinkers, in order to prove Theorem \ref{Gauss-space-2} we need some notion of weak holomorphy which can can be used for self-shrinkers. This is given by the following result, which generalizes Theorem \ref{E-T-0}, also proved by the authors in \cite{ANZ-1}:

\begin{theorem}[Theorem 1.1 of \cite{ANZ-1}]\label{main}
Let $h:U\subset\CC\rightarrow\CC$ be a complex function defined in an open set $U$ of the complex plane and $z=z_0\in U$ be a zero of $h.$ If there exists $\vp\in L^p_{loc}(U),\ p>2,$ a non-negative real function such that 
\begin{equation}\label{cauchy}
\left|\frac{\partial h}{\partial\bar{z}}\right|\leq \vp(z)G(|h(z)|),
\end{equation}
where $G:[0,\infty)\rightarrow[0,\infty)$  is a locally integrable function such that \[\limsup_{t\rightarrow0^+}\frac{G(t)}{t}<\infty,\] then either $h=0$ in a neighborhood $V\subset U$ of $z_0,$ or
\begin{equation}\label{weak-h}
h(z)=(z-z_0)^k h_k(z), \ z\in V, \ k\geq1,
\end{equation}
where $h_k(z)$ is a continuous function with $h_k(z_0)\neq 0.$
\end{theorem}

%\begin{remark}
%The existence of a weak notion of holomorphy to conclude \eqref{weak-h} was noticed first, as we know, by Carleman in 1933, see \cite{Carleman}, and was used later by Hartman and Wintner \cite{H-W} and \cite{H-W-2}, Chern \cite{chern}, Eschenburg and Tribuzy \cite{E-T-1} and \cite{E-T}, and Alencar-do Carmo-Tribuzy \cite{AdCT}.
%\end{remark}

Theorem \ref{main} has the following immediate consequence:

\begin{corollary}
Let $h:U\subset\CC\rightarrow\CC$ be a complex function defined in an open set $U$ of the complex plane. If \eqref{cauchy} holds, then on each connected components of $U$ which contains a zero of $h$, either $h\equiv 0$ or the  zeroes of $h$ are isolated.
\end{corollary}

\begin{remark}
%{\normalfont
The case when $\vp=0$ is equivalent to that $h$ is holomorphic. The case when $G(t)=t$ and $\vp$ is continuous, Theorem \ref{main} is the Main Lemma in \cite{AdCT} which implies Chern's Lemma in \cite{chern}.  Theorem \ref{main} also implies  Lemma 2.3, p. 154, of \cite{E-T}. There are many functions satisfying the condition $\limsup_{t\rightarrow0}G(t)/t<\infty.$ In fact, if $G$ is a continuous function such that $G(0)=0$, then $\limsup_{t\rightarrow0}G(t)/t=G'(0),$ if it exists. Moreover, if $G$ is any convex function with $G(0)=0$, then $G(t)/t\leq G(1)$ for small $0<t<1,$ which implies that convex functions also satisfy the condition. In particular, the functions $G(t)=t^\alpha, \alpha\geq 1,$ satisfy the condition. On the other hand, there are concave functions which satisfy this condition, for example $G(t)=\sin t,\ 0\leq t \leq \pi/2.$
\end{remark}

Theorem \ref{Gauss-space-2} is a consequence of the more general result Theorem \ref{Gauss-space}, p.\pageref{Gauss-space}, which holds for parallel weighted mean curvature surfaces in $\R^{2+m},\ m\geq1,$ where the weight is a radial function (i.e., which depends only on the distance of the point to the origin). In order to state this result, we shall need to give a brief introduction to weighted geometry in $\R^n.$ We refer, for example, \cite{C-M-Z} for a more detailed exposition. We call ($\R^n,\lan\cdot,\cdot\ran,e^{-f}$) a weighted Riemannian manifold if it has a weighted measure $dV_f = e^{-f}dV,$ where $f:\R^n\rightarrow\R$ is a function of class $C^2$. Let $X:\Sigma\rightarrow\R^n$ be an immersion of a surface $\Sigma.$ Consider $\Sigma$ with the weighted measure 
\[
 d\Sigma_f = e^{-f}d\Sigma,
\]
and the induced metric $\lan\cdot,\cdot\ran.$ 

The first variation of the weighted volume $V_f(\Sigma)=\int_\Sigma e^{-f}d\Sigma$ is given by
\[
\left.\dfrac{d}{dt}V_f(\Sigma_t)\right|_{t=0} = -\int_\Sigma \lan T^{\perp}, {\bf H}_f\ran e^{-f}d\Sigma,
\]
where $T$ is a compactly supported variational vector field on $\Sigma$ and
\begin{equation}\label{H_F}
{\bf H}_f = {\bf H} + (\n f)^\perp 
\end{equation}
is the weighted mean curvature vector of $\Sigma$ in $\R^n.$ Here $(\n f)^\perp$ denotes the part of the gradient $\n f $ of $f$ in $\R^n$ normal to $\Sigma$ and ${\bf H}$ denotes the non-normalized mean curvature vector of $\Sigma$ in $\R^n,$ i.e., the trace of the operator
\[
B(Z,W) = \n_ZW - \n^\Sigma_ZW,
\]
where $\n$ and $\n^\Sigma$ denote the connections of $\R^n$ and $\Sigma,$ respectively.

We say that a surface $\Sigma$ has parallel weighted mean curvature if ${\bf H}_f$ is parallel in the normal bundle, i.e., $\n^\perp{\bf H}_f=0.$ In particular, if ${\bf H}_f=0,$ we say that $\Sigma$ is $f$-minimal.

In the case that $f(X)=\|X\|^2/4,$ we call the weighted manifold \\($\R^n,\lan\cdot,\cdot\ran,e^{-\|X\|^2/4}$) the Gaussian space. Notice that self-shrinkers are $f$-minimal surfaces in the Gaussian space.

If the codimension is one, the parallel weighted mean curvature surfaces in the Gaussian space are called $\lambda$-surfaces. By using \eqref{H_F}, we can see that $\lambda$-surfaces are characterized by the equation
\[
\lambda = H + \frac{1}{2}\lan X,N\ran,
\]
where $\lambda\in\R,$ $N$ is the unit normal vector field of the immersion, and $H$ is its mean curvature, i.e., ${\bf H}=HN.$ Observe that self-shrinkers of $\R^3$ are also $\lambda$-surfaces for $\lambda=0.$

For each point $p\in\Sigma$ we can take isothermal parameters $u$ and $v$ in a neighborhood of $p,$ i.e., \[ds^2=\alpha(u,v)(du^2+dv^2),\] where $ds^2$ is the metric of $\Sigma$ and $\alpha(u,v)$ is a positive smooth function on $\Sigma$. Complexifying the parameters by taking $z=u+iv,$ we can identify $\Sigma$ with a subset of $\mathbb{C}.$ In this case, we have
\[
\lan X_z,X_{\bar{z}}\ran = \frac{\alpha(z)}{2} \ \mbox{and} \ ds^2=\alpha(z)|dz|^2.
\]
The immersion $X$ satisfies the equations
\begin{equation}\label{codazzi}
\left\{
\begin{aligned}
\n_{X_z}X_z&=\frac{\alpha_z}{\alpha}X_z + B(X_z,X_z),\\
\n_{X_{\bar{z}}}X_z&=\frac{\alpha}{4}{\bf H},\\
\n_{X_{\bar{z}}}X_{\bar{z}}&= \frac{\alpha_{\bar{z}}}{\alpha}X_{\bar{z}} + B(X_{\bar{z}},X_{\bar{z}}),\\
\end{aligned}
\right.
\end{equation}
and, for any $\nu\in T\Sigma^\perp,$
\begin{equation}\label{codazzi-2}
\left\{
\begin{aligned}
\n_{X_z}\nu &=-\frac{1}{2}\lan{\bf H},\nu\ran X_z - \frac{2}{\alpha}\lan B(X_z,X_z),\nu\ran X_{\bar{z}} + \n^\perp_{X_z}\nu\\
\n_{X_{\bar{z}}}\nu&= - \frac{2}{\alpha}\lan B(X_{\bar{z}},X_{\bar{z}}),\nu\ran X_z - \frac{1}{2}\lan{\bf H},\nu\ran X_{\bar{z}} + \n^\perp_{X_{\bar{z}}}\nu,\\
\end{aligned}
\right.
\end{equation}
where $\n^\perp$ is the connection of the normal bundle $T\Sigma^\perp.$

Let us denote by
\[
P^\nu dz^2 = \lan B(X_z,X_z),\nu\ran dz^2
\]
the $(2,0)$-part of the second fundamental form of $\Sigma$ in $\R^n$ relative to the normal $\nu\in T\Sigma^\perp.$ This quadratic form is also called the Hopf quadratic differential.

The follolwing theorem, which was proven in \cite{ANZ-1}, is a rigidity result for parallel weighted mean curvature ${\bf H}_f$ surfaces in the Euclidean space with arbitrary codimension and radial weight $f(X)=F(\|X\|^2),$ where $F:\R\to\R$ is a function of class $C^2$. Since the codimension can be  arbitrary large, we  assume that $X(\Sigma)$ does not lie in any proper affine subspace of the Euclidean space.

\begin{theorem}\label{Gauss-space}
Let $X:\Sigma\to\R^{2+m},\ m\geq1,$ be an immersion of a surface homeomorphic to the sphere. Assume that all the following assertions holds:
\begin{itemize}
\item[i)] $X$ has parallel weighted mean curvature ${\bf H}_f,$ i.e., $\n^\perp{\bf H}_f=0,$ for a radial weight $f(X)=F(\|X\|^2)$, where $F:\R\to\R$ is a function of class $C^2$.
\item[ii)] There exists a unitary normal vector field $\nu\in T\Sigma^\perp$ such that $\n^\perp\nu=0.$
\item[iii)] There exists a non-negative locally $L^p$ function $\vp:\Sigma\to\R,\ p>2,$ and a locally integrable function $G:[0,\infty)\rightarrow[0,\infty)$ satisfying $\limsup_{t\rightarrow 0}G(t)/t<\infty,$ such that  
\begin{equation}\label{HG}
\begin{aligned}
\left|F'(\|X\|^2)\lan {\bf H}_f,\nu\ran - 2\left[2F''(\|X\|^2)+(F'(\|X\|^2))^2\right]\right.&\left.\lan X,\nu\ran\right|\|X^\top\|\\
&\leq \vp G(\|\Phi^\nu\|).
\end{aligned}
\end{equation}
\end{itemize}
Then $X(\Sigma)$ is contained in a round hypersphere of $\R^{2+m}.$ Moreover, if ${\bf H}\neq 0$ and $\nu={\bf H}/\|{\bf H}\|,$ then $X(\Sigma)$ is a minimal surface of a round hypersphere of $\R^{2+m}$ or it is a round sphere in $\R^{2+m}.$ 

Here $X^\top$ denotes the component of $X$ tangent to $T\Sigma,$ $\|\Phi^\nu\|$ denotes the matrix norm of $\Phi^\nu=A^\nu-(\tr A^\nu/2)I,$ where $A^\nu$ is the shape operator of the second fundamental form of $X$ relative to $\nu,$ $\tr A^\nu$ is its trace, and $I:T\Sigma\to T\Sigma$ is the identity operator.
\end{theorem}

The ideia of the proof is to apply the Cauchy-Riemann inequality of Theorem \ref{main}, to the quadratic differential $Q^\nu = e^{-\frac{1}{2}f}P^\nu$ and conclude that $P^\nu,$ is identically zero in a neighborhood $V$ of a zero $z_0$ or this zero is isolated and the index of a direction field determined by $\textrm{Im}[P^\nu dz^2]=0$ is negative. If, for some coordinate neighborhood $V$ of zero, $P^\nu=0$, this holds for the whole $\Sigma$. Otherwise, the zeroes on the boundary of $V$ will contradict to Theorem \ref{main}. So if $P^\nu$ is not identically zero, all zeroes, if any, are isolated and have negative indices. This implies that the sum of all indexes of the isolated zeroes are negative (if there are zeroes) or zero (if there are no zeroes). Since $\Sigma$ has genus zero, by the Poincar\'e index theorem the sum of the indices of the singularities of any field of directions is $2$ (hence positive). This contradiction shows that $P^\nu$ is identically zero. This implies that $A^\nu = \mu I,$ i.e., $\nu$ is a umbilical normal direction of $X$. We then prove that $\mu$ must be constant and, since $X(\Sigma)$ does not lies in a hyperplane, we conclude that $\mu\neq 0$ and $X(\Sigma)$ lies in a hypersphere of $\R^{2+m}.$ This fact comes from Yau \cite{Yau} (see Theorem 1, p.351-352) and Chen-Yano\cite{Chen-Yano} (see Theorem 3.3, p.472-473). Moreover, if ${\bf H}$ and $\nu={\bf H}/\|{\bf H}\|,$ then $X(\Sigma)$ is a minimal surface of a hypersphere of $\R^{2+m}.$ This comes from Theorem 2, p.117, of the work of Ferus \cite{Ferus}.

In the case when $\Sigma$ is $f$-minimal, i.e., ${\bf H}_f=0,$ and the weight $f(X)=F(\|X\|^2)$ satisfies $F'(t)\neq 0$ and $2F''(t) + (F'(t))^2\neq 0,$ for every $t\in\R,\ t\geq0,$ the next result follows from Theorem \ref{Gauss-space}.

\begin{corollary}\label{Gauss-space-cor}
Let $X:\Sigma\to\R^{2+m},\ m\geq1,$ be an immersion of a surface homeomorphic to the sphere. Assume that all the following assertions holds:
\begin{itemize}
\item[i)] $X$ is $f$-minimal, i.e., ${\bf H}_f=0,$ for a radial weight $f(X)=F(\|X\|^2)$, where $F:\R\to\R$ is a function of class $C^2$ such that $F'(t)\neq0$ and $2F''(t) + (F'(t))^2\neq 0,$ for every $t\in\R, \ t\geq 0.$
\item[ii)] There exists an unitary normal vector field $\nu\in T\Sigma^\perp$ such that $\n^\perp\nu=0.$
\item[iii)] There exists a non-negative locally $L^p$ function $\vp:\Sigma\to\R,\ p>2,$  and  a locally integrable function $G:[0,\infty)\rightarrow[0,\infty)$ satisfying $\limsup_{t\rightarrow 0}G(t)/t<\infty,$ such that  
\begin{equation}\label{HG-cor}
\left(\|X\|^2 - \left(\frac{\|{\bf H}\|}{2F'(\|X\|^2)}\right)^2\right)\left(\frac{|\lan{\bf H},\nu\ran|}{2F'(\|X\|^2)}\right)^2\leq \vp^2 G(\|\Phi^\nu\|)^2.
\end{equation}
\end{itemize}

Then $X(\Sigma)$ is contained in a round hypersphere of $\R^{2+m}$ of radius $R,$ where $R$ is the solution of the equation 
\[
F'(R^2)R^2=1,
\]
and centered at the origin. Moreover, if ${\bf H}\neq 0$ and $\nu={\bf H}/\|{\bf H}\|,$ then $X(\Sigma)$ is a minimal surface of a round hypersphere of $\R^{2+m}$ with the same properties.

Here $\|\Phi^\nu\|$ is the matrix norm of $\Phi^\nu=A^\nu-(\tr A^\nu/2)I,$  where $A^\nu$ is the shape operator of the second fundamental form of $X$ relative to $\nu,$ $\tr A^\nu$ is its trace, and $I:T\Sigma\to T\Sigma$ is the identity operator.
\end{corollary}

Since self-shrinkers are $f$-minimal surfaces for the weight $f(X)=\frac{\|X\|^2}{4},$ applying Corollary \ref{Gauss-space-cor} to $F(t)=t/4,$ we obtain

\begin{corollary}\label{codim}
Let $X:\Sigma\rightarrow\R^{2+m},\ m\geq1,$ be an immersed self-shrinker homeomorphic to the sphere. Assume there exists an unitary normal vector field $\nu\in T\Sigma^\perp$ such that $\n^\perp\nu=0.$ If there exists a non-negative locally $L^p$ function $\vp:\Sigma\to\R,\ p>2,$ and a locally integrable function $G:[0,\infty)\rightarrow[0,\infty)$ satisfying $\limsup_{t\rightarrow 0}G(t)/t<\infty,$ such that  
\[
\left(\|X\|^2 - 4\|{\bf H}\|^2\right)|\lan{\bf H},\nu\ran|^2\leq \vp^2 G(\|\Phi^\nu\|)^2,
\]
then $X(\Sigma)$ is contained in a round hypersphere of $\R^{2+m}$ of radius $2$ and centered at the origin.

Here $\|\Phi^\nu\|$ is the matrix norm of $\Phi^\nu=A^\nu-(\tr A^\nu/2)I,$ where $A^\nu$ is the shape operator of the second fundamental form of $X$ relative to $\nu,$ $\tr A^\nu$ is its trace, and $I:T\Sigma\to T\Sigma$ is the identity operator.
\end{corollary}

If we consider the case of codimension one in Corollary \ref{codim}, then we obtain Theorem \ref{Gauss-space-2}:

\begin{corollary}[Theorem \ref{Gauss-space-2}]
Let $X:\Sigma\to\R^3$ be an immersed self-shrinker homeomorphic to the sphere. If there exists a non-negative locally $L^p$ function $\vp:\Sigma\to\R,\ p>2,$ and a locally integrable function $G:[0,\infty)\rightarrow[0,\infty)$ satisfying $\limsup_{t\rightarrow 0}G(t)/t<\infty,$ such that  
\[
(\|X\|^2-4H^2)H^2\leq \vp^2 G(\|\Phi\|)^2,
\]
then $X(\Sigma)$ is a round sphere of radius $2$ and centered at the origin.

Here  $\|\Phi\|$ denotes the matrix norm of $\Phi=A-(H/2)I,$ where $A$ is the shape operator of the second fundamental form of $X,$ $H$ is its non-normalized mean curvature, and $I$ is the identity operator of $T\Sigma.$
\end{corollary}

For surfaces with parallel weighted mean curvature in the Gaussian space, we have 

\begin{corollary}\label{Cor-lambda-000}
Let $X:\Sigma\rightarrow(\mathbb{R}^{2+m},\lan\cdot,\cdot\ran,e^{-\|X\|^2/4}),\ m\geq1,$ be an immersion of a surface homeomorphic to the sphere into the Gaussian space. Assume that all the following assertions holds:
\begin{itemize}
\item[i)] $X$ has parallel weighted mean curvature ${\bf H}_f,$ i.e., $\n^\perp{\bf H}_f=0$.
\item[ii)] There exists an unitary normal vector field $\nu\in T\Sigma^\perp$ such that $\n^\perp\nu=0.$
\item[iii)] There exists a non-negative locally $L^p$ function $\vp:\Sigma\to\R,\ p>2,$ and a locally integrable function $G:[0,\infty)\rightarrow[0,\infty)$ satisfying $\limsup_{t\rightarrow 0}G(t)/t<\infty,$ such that  
\begin{equation}\label{Cor-lambda-0}
(\|X\|^2 -4\|{\bf H}_f - {\bf H}\|^2)\lan {\bf H},\nu\ran^2\leq\varphi^2 G(\|\Phi^\nu\|)^2,
\end{equation}
\end{itemize}

Then $X(\Sigma)$ is contained in a round hypersphere of $\R^{2+m}.$ Moreover, if ${\bf H}\neq0$ and $\nu={\bf H}/\|{\bf H}\|,$ then $X(\Sigma)$ is a minimal surface of a round hypersphere of $\R^{2+m}$ of radius 
\[
\sqrt{\lan{\bf H}_f,\nu\ran^2 + 4} - \lan{\bf H}_f,\nu\ran.
\]

Here $\|\Phi^\nu\|$ is the matrix norm of $\Phi^\nu=A^\nu-(\tr A^\nu/2)I,$ where $A^\nu$ is the shape operator of the second fundamental form of $X$ relative to $\nu,$ $\tr A^\nu$ is its trace, and $I:T\Sigma\to T\Sigma$ is the identity operator.
\end{corollary}

In particular, for $\lambda$-surfaces, we obtain
\begin{corollary}

Let $X:\Sigma\rightarrow\R^3$ be a immersed $\lambda$-surface homeomorphic to the sphere. If there exists a non-negative locally $L^p$ function $\vp:\Sigma\to\R,\ p>2,$ and a locally integrable function $G:[0,\infty)\rightarrow[0,\infty)$ satisfying $\limsup_{t\rightarrow 0}G(t)/t<\infty,$ such that  
\[
\left(\|X\|^2-4(\lambda-H)^2\right)H^2\leq\varphi^2 G(\|\Phi\|)^2,
\]
then $X(\Sigma)$ is a round sphere of radius $\sqrt{\lambda^2+4}-\lambda$ and center at the origin.

Here  $\|\Phi\|$ denotes the matrix norm of $\Phi=A-(H/2)I,$ where $A$ is the shape operator of the second fundamental form of $X,$ $H$ is its non-normalized mean curvature, and $I$ is the identity operator of $T\Sigma.$
\end{corollary}

\begin{remark}
In the proof of Corollary \ref{Cor-lambda-000}, since the codimension can be $m\geq 2,$ we have that the spheres $\|X\|^2=constant$ and $\mathbb{S}^{1+m}(x_0,R)$ could be different. In this case we will have 
\[
X(\Sigma)\subset \mathbb{S}^{1+m}(x_0,R)\cap\mathbb{S}^{1+m}(0,\|X\|),
\]
where this intersection is, by its turn, a $m$-dimensional sphere.
\end{remark}

\section{Self-shrinkers of extrinsic curvature flows}

In this section, we continue the work done in \cite{ANZ-1} presenting results in the same spirit of Theorem \ref{Gauss-space-2} to other famous curvature flows, as the flow by the powers of Gaussian curvature $K$, the flow by the powers of the harmonic mean curvature $K/H$ and the flow by the powers of the mean curvature $H$. All these flows can be seen as particular cases of the general curvature flow, as follows.

Given a initial immersion $X_0:\Sigma\to\R^3$ of a two-dimensional surface, we say that the evolution of $X_0(\Sigma)$ by the curvature is a smooth one-parameter family of immersions $X:\Sigma\times[0,T)\to\R^3$ satisfying the initial value problem \begin{equation}\label{F-1}
\begin{cases}
\dfrac{\partial X}{\partial t} = W(k_1,k_2)N,\\
X(\cdot,0)=X_0,\\
\end{cases}
\end{equation}
where $k_1$ and $k_2$ are the principal curvatures of the immersions $X$, $N$ is their unitary normal vector fields, and $W\in C^1(\R^2)$. It is known this flow will be parabolic if and only if
\begin{equation}\label{pl}
\frac{\partial W}{\partial k_1}\cdot \frac{\partial W}{\partial k_2} >0.
\end{equation}

The more important curvature flows are those whose function $W$ is a combination of the mean curvature $H=k_1+k_2$ and the Gaussian curvature $K=k_1k_2.$ Among these we can cite the mean curvature flow, for $W(k_1,k_2)=H,$ the Gaussian curvature flow $W(k_1,k_2)=K,$ and the harmonic curvature flow, for $W(k_1,k_2)=K/H.$ These curvature flows has been studied by many authors in the last three decades, see \cite{Gerhardt-1}, \cite{Urbas-0}, \cite{Urbas-01}, \cite{Andrews-FN-0}, \cite{Andrews-FN-01}, \cite{Andrews-FN-1}, \cite{Andrews-FN-2}, \cite{McCoy} and references therein.

In the study the curvature flows, the self-similar solutions play an important role since its was proved that, under some convexity conditions, the solutions of the flow, when suitably normalized, converge to a self-similar solution. A solution of the curvature flow \eqref{F-1} is said self-similar if each $\Sigma_t=X(\Sigma,t)$ is an homothety, or a translation, or even a rotation of $\Sigma_0=X_0(\Sigma).$ The homothetic self-similar solutions are said self-shrinkers (or shrinking self-similar solutions, or shrinking homothetic solutions) if the solution shrinks homothetically from $\Sigma_0.$ 

If $W$ is a homogeneous function of degree $\beta>0,$ i.e., $W(ak_1,ak_2)=a^\beta W(k_1,k_2),\ a>0,$ then it can be proved that a shrinking self-similar solution of a curvature flow satisfies the equation
\begin{equation}\label{soliton}
W(k_1,k_2)=-\lambda\lan X,N\ran, \ \lambda\in(0,\infty).
\end{equation}
Changing the variables $x_1=k_1+k_2$ and $x_2=(k_1-k_2)^2$ we can write
\begin{equation}\label{W-Psi0}
W(k_1,k_2)=\Psi(x_1,x_2)=\Psi(k_1+k_2,(k_1-k_2)^2)=\Psi(H,H^2-4K),
\end{equation}
where $K=k_1k_2$ is the Gaussian curvature of the immersion $X$ and $H=k_1+k_2$ is its mean curvature. Therefore, the equation \eqref{soliton} becomes
\begin{equation}\label{soliton-0}
\Psi(H,H^2-4K)=-\lambda\lan X,N\ran, \ \lambda\in(0,\infty).
\end{equation}

The main result in this subject is the following

\begin{theorem}[Alencar-Silva Neto-Zhou]\label{thm-soliton-0}
Let $X:\Sigma\to\R^3$ be a closed, immersed surface of genus zero satisfying \eqref{soliton-0}, where $\Psi: \mathbb{R}\times [0,+\infty)\to \mathbb{R}$ is a $C^1$ function satisfying $\frac{\partial\Psi}{\partial x_1}\neq 0$. If there exists a non-negative function $\vp\in L^p(\Sigma),$ $p>2,$ and a locally integrable function $G:[0,\infty)\rightarrow[0,\infty)$ satisfying $\limsup_{t\rightarrow 0}G(t)/t<\infty,$ such that
\begin{equation}\label{hyp-soliton-0}
H^2(\|X\|^2-\lan X,N\ran^2)\leq \vp^2 G(\|\Phi\|)^2,
\end{equation}
then $X(\Sigma)$ is a round sphere centered at the origin and radius satisfying the equation 
\[
\lambda R = \Psi\left(\frac{2}{R},0\right).
\]
Here $\|\Phi\|$ denotes the matrix norm of $\Phi=A-(H/2)I,$ where $A$ is the shape operator of the second fundamental form of $X,$ $H$ is its non-normalized mean curvature, and $I$ is the identity operator of $T\Sigma.$
\end{theorem}

\begin{remark}
Theorem \ref{thm-soliton-0} cannot be derived from Theorem \ref{Gauss-space} of the previous section since the only curvature flow that can be expressed using weighted geometry is the mean curvature flow.
\end{remark}

\begin{remark}\label{rem-K1}
Since the function $\vp$ is assumed to be only $L^p,$ it is allowed to be infinity in some points. In particular, it is possible to have 
\[
\lim_{p\to p_0} \vp (p) G(\|\Phi (p)\|)>0
\]
for umbilical points $p_0\in\Sigma$ despite $G(\|\Phi (p_0)\|)=G(0)=0.$ Therefore, the inequality \eqref{hyp-soliton-0} does not imply necessarily that $H^2(\|X\|^2-\lan X,N\ran^2)=0$ at umbilical points.
\end{remark}

\begin{remark}\label{rem-HK}
Notice that, since $H^2-4K=(k_1-k_2)^2\geq 0,$ the inequality
\[
K\leq\frac{1}{4}H^2
\]
holds for every surface in $\R^3.$ As it was shown in the Remark 1.4 of \cite{ANZ-1}, inequality \eqref{hyp-soliton-0} gives the existence of a function $\psi,$ which can be chosen satisfying $\psi^2<\ve$ for every given $\ve>0$ arbitrarily small, such that $1/\psi\in L^p(\Sigma),$ $p>2,$ and
\[
K\leq\frac{1}{4}(1-\psi^2)H^2.
\]
\end{remark}

\begin{remark}
The hypothesis \eqref{hyp-soliton-0} of Theorem \ref{thm-soliton-0} is necessary. In fact, we prove in \cite{ANZ-2} that, if there exists non-spherical genus zero rotational surface which is the solution of \eqref{soliton-0}, then \eqref{hyp-soliton-0} does not hold.
\end{remark}

\begin{remark}
The function $W(k_1,k_2)$ is homogeneous of degree $\beta\in\R,$ if and only if the function $\Psi$ satisfies 
\begin{equation}\label{psi-homo}
\Psi(a x_1,a^2 x_2)=a^\beta\Psi(x_1,x_2),\ a>0.
\end{equation} 
By an abuse of notation, we will call $\Psi$ a homogeneous function of degree $\beta\in\R$ if $\Psi$ satisfies \eqref{psi-homo}. If $\Psi$ is homogeneous of degree $\beta\neq -1$ with $\Psi(1,0)>0,$ and $\lambda>0,$ then the radius of the sphere of Theorem \ref{thm-soliton-0} is given by
\[
R=\left[\lambda^{-1}2^\beta\Psi(1,0)\right]^{\frac{1}{\beta+1}}.
\]
\end{remark}

\begin{remark}
{\normalfont
The flow \eqref{F-1} is a (weakly) parabolic equation if and only if
\[
\frac{\partial W}{\partial k_1}\cdot\frac{\partial W}{\partial k_2}>0\ (\geq 0),
\]
or equivalently
\[
\left(\frac{\partial\Psi}{\partial x_1}\right)^2 - 4x_2\left(\frac{\partial\Psi}{\partial x_2}\right)^2>0\ (\geq 0).
\]
Notice the hypothesis $\frac{\partial\Psi}{\partial x_1}\neq 0$ of Theorem \ref{thm-soliton-0} assures the parabolicity of the flow near the umbilical points ($x_2=0$), but this result holds even when the flow is not parabolic.
}
\end{remark}

Our first consequence is for the $\alpha$-mean curvature flow, $\alpha\in \R\backslash\{0,1\},$
\[
\frac{\partial X}{\partial t}=H^\alpha N.
\]
This flow is parabolic for $H>0.$ The case when $\alpha=1$ is the so called mean curvature flow, which is parabolic without any additional assumption. This case was dealt by the authors in \cite{ANZ-1}. 

Schulze, see \cite{Scz-1}, proved that closed (weakly) convex hypersurfaces of $\R^{n+1}$ converges to a point if $\alpha\in(0,1)$ ($\alpha\geq 1$) and Schn\"urer, see \cite{SA2}, and Schulze, see \cite{Scz-2}, proved that closed convex surfaces of $\R^3$ converges to a round point for $1\leq \alpha \leq 5.$ For general speeds of higher homogeneity, Andrews, see \cite{Andrews-FN-2}, proved that the flow of a closed convex surface converges to a round point provided it satisfies an initial pinching condition. 

The shrinking self-similar solutions of the $\alpha$-mean curvature flow satisfy the equation
\[
H^\alpha = -\lambda\lan X,N\ran,\ \lambda\in (0,\infty).
\]
Our result characterizes the sphere as the only mean convex (i.e., $H\neq 0$), genus zero, closed shrinking self-similar solution of the $\alpha$-mean curvature flow under an upper pinching curvature condition. Notice that the mean convex assumption is weaker than convexity, since mean convexity admits immersed surfaces and surfaces with $K\leq 0.$
 
\begin{corollary}\label{Hn-alpha}
Let $X:\Sigma\to\R^3$ be a closed, homeomorphic to the sphere, immersed, mean convex, two-dimensional shrinking self-similar solution of the $\alpha$-mean curvature flow, for $\alpha\in \R\backslash\{-1,0,1\}.$ If there exists a non-negative function $\vp\in L^p(\Sigma),$ $p>2,$ and a locally integrable function $G:[0,\infty)\rightarrow[0,\infty)$ satisfying $\limsup_{t\rightarrow 0}G(t)/t<\infty,$ such that    
\begin{equation}\label{eq11}
H^2(\|X\|^2-\lan X,N\ran^2)\leq \vp^2 G(\|\Phi\|)^2,
\end{equation}
then $X(\Sigma)$ is a round sphere of radius $(2^\alpha\lambda^{-1})^{\frac{1}{\alpha+1}}$ and center at the origin.

Here $\|\Phi\|$ denotes the matrix norm of $\Phi=A-(H/2)I,$ where $A$ is the shape operator of the second fundamental form of $X,$ $H$ is its non-normalized mean curvature, and $I$ is the identity operator of $T\Sigma.$
\end{corollary}
\begin{remark}
{\normalfont
If $\alpha=\frac{m}{2n-1}\in(0,1),\ n,m\in\mathbb{N},$ then the hypothesis of mean convexity in Corollary \ref{Hn-alpha} is not necessary. Notice that in this case the flow is only weakly parabolic, becoming degenerate for the points when $H=0.$
}
\end{remark}
\begin{remark}
Drugan, Lee and Wheeler \cite{Drugan} proved that the spheres are the only closed self-shrinkers for the inverse mean curvature flow (i.e., for $\alpha=-1$) without any additional assumption, solving the problem in this case.
\end{remark}

The next application of Theorem \ref{thm-soliton-0} is for the $\alpha$-harmonic mean curvature flow
\[
\frac{\partial X}{\partial t} =\left(\frac{K}{H}\right)^\alpha N,
\]
whose shrinking self-similar solitons satisfy the equation
\[
\left(\frac{K}{H}\right)^\alpha=-\lambda\lan X,N\ran,\ \lambda\in (0,\infty).
\]
If $\alpha\in (0,\infty),$ then this flow is (weakly) parabolic for (weakly) convex surfaces, being degenerate for the points where $K=0.$ If we consider values of $\alpha$ such that $K$ can assume negative values, as $\alpha=\frac{m}{2n-1},\ m,n\in\mathbb{N},$ including the classical case of $\alpha=1,$ then the flow is weakly parabolic for every surface, being degenerate for the points where $K=0$ and singular for the points where $H=0.$

For $\alpha=1,$ the existence of solutions for closed convex surfaces as initial data was proved by Andrews, see \cite{Andrews-FN-0}, who also showed that closed convex surfaces flowing by the harmonic mean curvature converges to a round point in finite time (in fact, the result of Andrews holds for a more wide class of degree one homogeneous functions $W$). Dieter, see \cite{Dieter}, studied the convergence of the flow for the degenerate case $K\geq 0$ and $H>0,$ Caputo and Daskalopoulos, see \cite{D-C}, and Daskalopoulos and Sesum, see \cite{D-S}, studied the highly degenerate case, where $K$ and $H$ can be simultaneously zero. The case when $K<0$ and $H<0$ was studied by Daskalopoulos and Hamilton, see \cite{D-H-2}. 

For $\alpha\in(0,1),$ Anada, see \cite{Anada}, proved the existence of non-round closed convex self-similar solutions of the $\alpha$-harmonic mean curvature flow. After this findings, in a joint work with Tsutsumi, see \cite{A-T}, he also investigated sufficient conditions for the $\alpha$-mean curvature flows converge to a round point. 

Our result gives conditions for a closed, mean convex, self-similar solution with genus zero of the $\frac{m}{2n-1}$-harmonic mean curvature flow to be a sphere. We remark here that the powers $\frac{m}{2n-1}, m,n\in\mathbb{N},$ allows us to work with surfaces such that $K<0$ at some points, but our technique holds for every $\alpha\in(0,1],$ if we assume that $\Sigma$ is weakly convex.

\begin{corollary}\label{harm-H}
Let $X:\Sigma\to\R^3$ be a closed, homeomorphic to the sphere, immersed, mean convex, two-dimensional shrinking self-similar solution of the $\alpha$-harmonic mean curvature flow for $\alpha=\frac{m}{2n-1},$ where $m, n\in\mathbb{N}$ and $\frac{m}{2n-1}\leq 1.$ If there exists a non-negative function $\vp\in L^p(\Sigma),$ $p>2,$ and a locally integrable function $G:[0,\infty)\rightarrow[0,\infty)$ satisfying $\limsup_{t\rightarrow 0}G(t)/t<\infty,$ such that    
\begin{equation}\label{eq10}
H^2(\|X\|^2-\lan X,N\ran^2)\leq \vp^2 G(\|\Phi\|)^2,
\end{equation}
then $X(\Sigma)$ is a round sphere of radius $(2^{\frac{m}{2n-1}}\lambda)^{-\frac{2n-1}{m-2n+1}},$ centered at the origin, if $(m,n)\neq(1,1),$ and for any radius $R>0,$ centered at the origin, with $\lambda=\frac{1}{2},$ if $(m,n)=(1,1).$

Here $\|\Phi\|$ denotes the matrix norm of $\Phi=A-(H/2)I,$ where $A$ is the shape operator of the second fundamental form of $X,$ $H$ is its non-normalized mean curvature, and $I$ is the identity operator of $T\Sigma.$
\end{corollary}

The last classical flow we will discuss here and obtain consequences of Theorem \ref{thm-soliton-0} is the $\alpha$-Gaussian curvature flow 
\[
\frac{\partial X}{\partial t} = K^\alpha N,
\]
whose shrinking self-similar solutions satisfy the equation 
\[
K^\alpha = -\lambda\lan X,N\ran, \ \lambda\in(0,\infty).
\]
This flow is (weakly) parabolic if $K>0$ ($K\geq 0$), being degenerate for the points where $K=0.$

When $\alpha=1,$ this flow is called Gaussian curvature flow, and was first introduced by Firey in 1974, see \cite{Firey}, as a model of the wearing process of convex rolling stones on a beach. He proved also that closed convex surfaces under this flow converges to a round point when the initial surface is symmetric about the origin. Tso, see \cite{Tso}, for $\alpha=1,$ and Chow, see \cite{Chow-1}, for $\alpha=1/n,$ proved the convergence to a point of a closed convex hypersurfaces of $\R^{n+1}$ under the flow. Andrews, see \cite{Andrews-1996}, proved that, for $\alpha=1/(n+2),$ closed convex hypersurfaces evolving under the flow converges to an ellipsoid. We observe that Calabi, see \cite{Calabi}, early proved that the ellipsoids are the only closed hypersurfaces satisfying the equation of the self-similar solutions of the $\frac{1}{n+2}$-Gaussian curvature flow. The works of Andrews, see \cite{Andrews-2000}, \cite{AGN}, and Guan and Ni, see \cite{Guan-Ni}, proved that the flow converges to a self-similar solution for every $\alpha\geq 1/(n+2).$ To conclude the analysis of the case when $\alpha\geq 1/(n+2),$ Brendle, Choi and Daskalopoulos, see \cite{B-C-D}, proved that the only closed self-similar solutions of the $\alpha$-Gaussian curvature flow for $\alpha>1/(n+2)$ are the round spheres. In his turn, if $\alpha<0,$ then Gerhardt, see \cite{Gerhardt-2}, proved that the only closed convex self-similar solution of the $\alpha$-Gaussian curvature flow is a round sphere. Moreover, he proved that the flow converges to a sphere after rescaling. 

On the other hand, Andrews, see \cite{Andrews-2000}, proved the existence of non-spherical closed convex self-similar solutions of the $\alpha$-Gaussian curvature flow for small $\alpha>0.$ In particular, in dimension $2,$ for $\alpha\in (0,1/10).$ This shows that if we want to characterize the sphere as the only self-similar solution of the $\alpha$-Gauss curvature flow for small values of $\alpha>0,$ then we will need some additional assumption. 

Our result provides sufficient conditions for a self-similar solution of the $\alpha$-Gaussian curvature flow, $\alpha\in(0,1/4),$ to be a round sphere.

\begin{corollary}\label{K-alpha}
Let $X:\Sigma\to\R^3$ be a closed, convex, two-dimensional shrinking self-similar solution of the $\alpha-$Gaussian curvature flow for $\alpha\in(0,1/4).$ If there exists a non-negative function $\vp\in L^p(\Sigma),$ $p>2,$ and a locally integrable function $G:[0,\infty)\rightarrow[0,\infty)$ satisfying $\limsup_{t\rightarrow 0}G(t)/t<\infty,$ such that  
\begin{equation}\label{eq9}
H^2(\|X\|^2-\lan X,N\ran^2)\leq \vp^2 G(\|\Phi\|)^2
\end{equation}
then $X(\Sigma)$ is a round sphere of radius $\lambda^{-\frac{1}{2\alpha+1}}$ and center at the origin.

Here $\|\Phi\|$ denotes the matrix norm of $\Phi=A-(H/2)I,$ where $A$ is the shape operator of the second fundamental form of $X,$ $H$ is its non-normalized mean curvature, and $I$ is the identity operator of $T\Sigma.$
\end{corollary}

\begin{remark}
{\normalfont
Since there are examples of closed convex self-similar solutions of the $\alpha$-Gaussian curvature flow for $\alpha\in(0,1/10),$ given by Andrews, see \cite{Andrews-2000}, at least in this cases some additional hypothesis like \eqref{eq9} is necessary to obtain the conclusions of Corollary \ref{K-alpha}.
}
\end{remark}

\begin{remark}
{\normalfont
Corollary \ref{K-alpha} holds in a more general setting: if we choose values of $\alpha$ which allows negative values of $K,$ as for example $\alpha=\frac{m}{2n-1}\leq 1,$ $m,n\in\mathbb{N},$ then we can assume only that $\Sigma$ is a closed mean convex surface with genus zero to obtain the same conclusion, despite the flow is not parabolic in this case.
}
\end{remark}

In order to illustrate the scope of situations to which the Theorem \ref{thm-soliton-0} can be applied in the context of the curvature flows, we give here a list of examples of homogeneous functions $W(k_1,k_2)$ such that the flow \eqref{F-1} is parabolic including negative values of $K.$

\begin{itemize}
\item[(i)] $W(k_1,k_2)=aH^2+bK,$ $a,b\in\R.$ The flow is parabolic for \[K>-\dfrac{2a(2a+b)}{b^2}H^2.\] In this case, $\frac{\partial\Psi}{\partial x_1}\neq 0$ if and only if $H\neq 0,$ i.e., the surface is mean convex. 

As a particular situation, we have $W(k_1,k_2)=|A|^2=k_1^2+k_2^2,$ by taking $a=1$ and $b=-2.$ In this case, the flow is parabolic for $K>0,$ i.e., for convex surfaces. This flow was studied by Schn\"urer in \cite{SA2}.

\item[(ii)] $W(k_1,k_2)=aH^{2\alpha}+bK^\alpha,$ $a,b>0,$ $\alpha=\frac{m}{2n-1}\geq1, m,n\in\mathbb{N}.$ The flow is parabolic for \[4a^2+2ab\left(\dfrac{K}{H^2}\right)^{\alpha-1}+b^2\left(\dfrac{K}{H^2}\right)^{2\alpha-1}>0.\] In this case, $\frac{\partial\Psi}{\partial x_1}\neq 0$ if and only if $H\neq 0,$ i.e., the surface is mean convex.

\item[(iii)] $W(k_1,k_2)=H^{\frac23}+bK^{\frac13},$ $b\in(0,2^{5/3}).$ The flow is parabolic for $K\neq 0$ and $H\neq 0.$ In this case, $\frac{\partial\Psi}{\partial x_1}\neq 0$ everywhere and it is singular for $K=0$ and $H=0.$
\end{itemize}

\section{Warped product manifolds}

In this section we present some results proved by the first and the second authors in \cite{ASN} which generalize the Eschenburg-Tribuzy theorem for the more general class of three-dimensional Riemannian manifolds $M^3=I\times\s^2,$ where $I=(0,b)$ or $I=(0,\infty),$ with the metric 
\begin{equation}\label{warped}
\lan\cdot,\cdot\ran =dt^2+h(t)^2d\om^2,    
\end{equation}
where $h:I\rightarrow\R$ is a smooth function, called warping function, and $d\om^2$ denotes the canonical metric of the $2$-dimensional round sphere $\s^2.$ With the metric \eqref{warped}, the product $M^3=I\times\s^2$ is called a warped product manifold and generalizes the space forms with constant sectional curvature. In fact, the metrics of the space forms of constant sectional curvature $c\in\R$ can be written in polar coordinates as in \eqref{warped}, where
\[
\begin{cases}
h(t)=t& \mbox{for} \ \R^3,\\
h(t)=\dfrac{1}{\sqrt{c}}\sin(\sqrt{c}t)& \mbox{for} \ \s^3(c),\\
h(t)=\dfrac{1}{\sqrt{-c}}\sinh(\sqrt{-c}t)& \mbox{for}\ \h^3(c).\\
\end{cases}
\]
The warped product manifold $M^3$ has two different sectional curvatures which depend only on the parameter $t$, one tangent to the slices $\{t\}\times\s^2,$ denoted by $K_{\tan}(t),$ and other relative to the planes which contains the radial direction $\partial t,$ which de denote by $K_{\rad}(t).$ In terms of the warping function, we can write
\begin{equation}\label{tan-rad}
K_{\tan}(t) =\dfrac{1-h'(t)^2}{h(t)^2} \ \mbox{and} \ K_{\rad}(t)=-\frac{h''(t)}{h(t)}.
\end{equation}

These manifolds were first introduced by Bishop and O' Neill in 1969, see \cite{B-ON}, and is having increasing importance due to its applications as model spaces in general relativity. Part of these applications comes from the metrics which are solutions of the Einstein equations, as the de Sitter-Schwarzschild metric and Reissner-Nordstrom metric, which we introduce later.

%In recent years, immersions in warped product manifolds have been extensively studied, with many interesting papers in this subject, for instance see \cite{montiel-1}, \cite{Montiel}, \cite{bray-2}, \cite{ritore}, \cite{B-M}, \cite{DR}, \cite{AD-1}, \cite{AD-2}, \cite{AIR}, \cite{bessa}, \cite{BCL}, \cite{Xia-Wu-1}, \cite{Gimeno}, \cite{GIR}, \cite{sal-sal}, \cite{Xia-Wu-2}, \cite{aledo}, \cite{GL}, \cite{ASN}, \cite{GLW}, and \cite{SN}. We can also cite the book of Petersen, see \cite{petersen}, for a modern presentation of warped product manifolds and the book of Besse \cite{besse} for an introduction to general relativity and the deduction of Schwarzschild space-time from the Einstein equations.

Applying the Hopf differential to Theorem \ref{CR-ET}, the main result of this section is the following generalization of Theorem \ref{E-T-0} for a class of warped product manifolds which contains the de Sitter-Schwarzschild and the Reissner-Nordstrom manifolds:

\begin{theorem}[Alencar-Silva Neto]\label{hopf-warped}
Let $\Sigma$ be a surface, homeomorphic to the sphere, immersed in a warped product manifold $M^3=I\times\s^2$, with mean curvature function $H$. If there exists a non-negative $L^p,$ $p>2,$ function $f:\Sigma\ria\R$ such that
\begin{equation}\label{ineq-warped}
\begin{aligned}
|dH + (K_{\tan}(t) &- K_{\rad}(t))\nu dt|\\
&\leq f\sqrt{H^2 - K + K_{\tan}(t) - (1-\nu^2)(K_{\tan}(t)-K_{\rad}(t))},
\end{aligned}
\end{equation}
then $\Sigma$ is totally umbilical. 

Moreover, if $K_{\tan}(t)\neq K_{\rad}(t),$ except possibly for a discrete set of values $t\in I,$ and $\Sigma$ has constant mean curvature, then $\Sigma$ is a slice.
\end{theorem}

\begin{remark}
Actually, some additional hypothesis as \eqref{ineq-warped} is needed in order to classify the slices as the only constant mean curvature spheres. In fact, it was observed by Brendle (see \cite{brendle-0}, Theorem 1.5, p. 250) that a result of Pacard and Xu (see \cite{pacard}, Theorem 1.1, p. 276) implies that in some warped product manifolds there are small spheres with constant mean curvature which are not umbilical.
\end{remark}

\begin{remark}
To obtain the slice in the second part of Theorem \ref{hopf-warped}, the assumption over $M^3$ that $K_{\tan}(t)\neq K_{\rad}(t),$ except possibly for a discrete set of values $t\in I,$ is necessary. In fact, if $K_{\tan}(t)=K_{\rad}(t)$ for some interval $(t_0,t_1)\subset I,$  then all the sectional curvatures of $M^3$ will depend only on $t.$ This will imply, by the classical Schur's Theorem, that $\widetilde{M}^3:=(t_0,t_1)\times\s^2$ has constant sectional curvature. In this case, there exists spheres, other than the slices, with constant mean curvature (in fact, the geodesic spheres centered in some point of $\widetilde{M}^3)$.
\end{remark}

Two of the most famous examples of warped product manifolds are the de Sitter-Schwarzschild manifolds and the Reissner-Nordstrom manifolds, which we describe below.

\begin{definition}[The de Sitter-Schwarzschild manifolds] Let $m>0,$ $c\in\R,$ and 
\[
(s_0,s_1)=\{r>0 ; 1-mr^{-1}-cr^2>0\}.
\] 
If $c\leq 0,$ then $s_1=\infty.$ If $c>0,$ assume that $cm^2<\frac{4}{27}.$ The de Sitter-Schwarzschild manifold is defined by $M^3(c)=(s_0,s_1)\times \s^2$ endowed with the metric
\[
\lan\cdot,\cdot\ran=\dfrac{1}{1-mr^{-1}-cr^2}dr^2 + r^2 d\om^2.
\]
In order to write the metric in the form \eqref{warped}, define $F:[s_0,s_1)\ria \R$ by
\[
F'(r)=\dfrac{1}{\sqrt{1-mr^{-1}-cr^2}}, \ F(s_0)=0.
\]
Taking $t=F(r),$ we can write $\lan\cdot,\cdot\ran=dt^2+h(t)^2d\om^2,$ where $h:[0,F(s_1))\ria[s_0,s_1)$ denotes the inverse function of $F.$ The function $h$ clearly satisfies
\begin{equation}\label{defi-SS}
h'(t)=\sqrt{1-mh(t)^{-1}-ch(t)^2},\ h(0)=s_0,\ \mbox{and}\ h'(0)=0.
\end{equation}
\end{definition}

For these manifolds, we have

\begin{corollary}[The de Sitter-Schwarzschild manifolds]
Let $\Sigma$ be a surface, homeomorphic to the sphere, immersed in the de Sitter-Schwarzschild manifold, with constant mean curvature. If there exists a non-negative $L^p,$ $p>2,$ function $f:\Sigma\ria\R$ such that
\[
|dt|\leq f \sqrt{H^2-K + c + \frac{m(3\nu^2-1)}{2h(t)^3}},
\]
then $\Sigma$ is a slice.

Here, $K$ is the Gaussian curvature of $\Sigma,$ $\nu=\lan \n t,N\ran$ is the angle function, and $N$ is the unitary normal vector field of $\Sigma$ in the de Sitter-Schwarzschild manifold.
% $\{t_0\}\times\s^2, t_0\in I.$
\end{corollary}

\begin{definition}[The Reissner-Nordstrom manifolds]\label{ex-RN} The Reissner- Nordstrom manifold is defined by $M^3=(s_0,\infty)\times\s^2,$ with the metric
\[
\lan\cdot,\cdot\ran=\dfrac{1}{1-mr^{-1}+q^2r^{-2}}dr^2 + r^2 d\om^2,
\]
where $m>2q>0$ and $s_0=\frac{2q^2}{m-\sqrt{m^2-4q^2}}$ is the larger of the two solutions of $1-mr^{-1}+q^2r^{-2}=0.$ In order to write the metric in the form \eqref{warped}, define $F:[s_0,\infty)\ria \R$ by
\[
F'(r)=\dfrac{1}{\sqrt{1-mr^{-1}+q^2r^{-2}}}, \ F(s_0)=0.
\]
Taking $t=F(r),$ we can write $\lan\cdot,\cdot\ran=dt^2+h(t)^2d\om^2,$ where $h:[0,\infty)\ria[s_0,\infty)$ denotes the inverse function of $F.$ The function $h$ clearly satisfies
\begin{equation}\label{defi-RN}
h'(t)=\sqrt{1-mh(t)^{-1}+q^2h(t)^{-2}},\ h(0)=s_0,\ \mbox{and}\ h'(0)=0.
\end{equation}
\end{definition}

For these manifolds, we have

\begin{corollary}[The Reissner-Nordstrom manifolds]
Let $\Sigma$ be a surface, homeomorphic to the sphere, immersed in the Reissner-Nordstrom manifold, with constant mean curvature. If there exists a non-negative $L^p,$ $p>2,$ function $f:\Sigma\ria\R$ such that
\[
|dt|\leq f \sqrt{H^2-K +\frac{m(3\nu^2-1)}{2h(t)^3}+ \frac{q^2(1-2\nu^2)}{h(t)^4}},
\]
then $\Sigma$ is a slice.

Here, $K$ is the Gaussian curvature of $\Sigma,$ $\nu=\lan \n t,N\ran$ is the angle function, and $N$ is the unitary normal vector field of $\Sigma$ in the Reissner-Nordstrom manifold.% $\{t_0\}\times\s^2, t_0\in I.$
\end{corollary}

\begin{remark}
Since the warped product manifold is smooth at $t=0$ if and only if $h(0)=0, \ h'(0)=1,$ and all the even order derivatives are zero at $t=0$, i.e., $h^{(2k)}(0) = 0,\ k > 0,$ see \cite{petersen}, Proposition 1, p. 13, we can see the de Sitter-Schwarzschild manifolds and the Reissner-Nordstrom manifolds are singular at $t=0.$
\end{remark}

\bibliographystyle{acm}
\bibliography{references}

\end{document}